\documentclass[12pt]{article}
\usepackage{amsmath,amsfonts}
\title{Visual Decompositions of Coxeter Groups}
\author{M. Mihalik\\and\\S. Tschantz}
\newtheorem{theorem}{Theorem}
\newtheorem{proposition}[theorem]{Proposition}
\newtheorem{lemma}[theorem]{Lemma}
\newtheorem{corollary}[theorem]{Corollary}
\newenvironment{proof}{\addvspace{12pt}\noindent{\bf Proof:}}{
$\Box$\par\addvspace{12pt}}
\newcounter{examplenum}
\newenvironment{example}{\addvspace{12pt}\refstepcounter{examplenum}
\noindent{\bf Example \arabic{examplenum}.}}{\par\addvspace{12pt}}
\newcounter{remarknum}
\newenvironment{remark}{\addvspace{12pt}\refstepcounter{remarknum}
\noindent{\bf Remark \arabic{remarknum}.}}{\par\addvspace{12pt}}
\date{December 20, 2005 VisualFinal2-05}
\begin{document}
\maketitle
\begin{abstract}
A Coxeter system is an ordered pair $(W,S)$ where $S$ is the
generating set in a particular type of presentation for the
Coxeter group $W$. A subgroup of $W$ is called special if it is
generated by a subset of $S$. Amalgamated product decompositions
of a Coxeter group having special factors and special amalgamated
subgroup are easily recognized from the presentation of the
Coxeter group.  If a Coxeter group is a subgroup of the
fundamental group of a given graph of groups, then the Coxeter
group is also the fundamental group of a graph of special
subgroups, where each vertex and edge group is a subgroup of a
conjugate of a vertex or edge group of the given graph of groups.
A vertex group of an arbitrary graph of groups decomposition of a
Coxeter group is shown to split into parts conjugate to special
groups and parts that are subgroups of edge groups of the given
decomposition. Several applications of the main theorem are
produced, including the classification of maximal FA-subgroups of
a finitely generated Coxeter group as all conjugates of certain
special subgroups.
\end{abstract}

\section{Introduction}\label{introsection}

We take a {\it Coxeter presentation} to be given as $$P=\langle
S:(st)^{m(s,t)}\ (s,t\in S,\,m(s,t)< \infty)\rangle$$ where
$m:S^2\to\lbrace 1,2,\ldots,\infty\rbrace$ is such that $m(s,t)=1$
iff $s=t$, and $m(s,t)=m(t,s)$.  In the group with this
presentation, the elements of $S$ represent distinct elements of
order 2 and a product $st$ of generators has order $m(s,t)$.  A
{\it Coxeter group} $W$ is a group having a Coxeter presentation
and a {\it Coxeter system} $(W,S)$ is a Coxeter group $W$ with
generating subset $S$ corresponding to the generators in a Coxeter
presentation of $W$. When the order of the product of a pair of
generators is infinite there will be no defining relator for that
pair of generators and we will say that the generators are {\it
unrelated}. A Coxeter group $W$ belongs to some Coxeter system
$(W,S)$, and though $S$ need not be uniquely determined up to an
automorphism of $W$, we often take such an $S$ as given. Our basic
reference for Coxeter groups is Bourbaki \cite{Bourbaki}. A {\it
special subgroup} of a Coxeter system $(W,S)$, is a subgroup of
$W$ generated by a subset of $S$ (see \cite{Brown}). If $W'$ is
the special subgroup generated by $S'\subseteq S$ in a Coxeter
system $(W,S)$, then $(W',S')$ is also a Coxeter system.

The information given by a Coxeter presentation may be
conveniently expressed in the form of a labeled graph. We define
the {\it presentation diagram} of the system $(W,S)$ to be the
labeled graph $\Gamma(W,S)$ with vertex set $S$, and an
(undirected) edge labeled $m(s,t)$ between distinct vertices $s$
and $t$ when $m(s,t)<\infty$. The connected components of the
presentation diagram $\Gamma(W,S)$ correspond to special subgroups
which are the factors in a free product decomposition of $W$.
(This is in contrast to a {\it Coxeter graph} with vertex set $S$
and labeled edges when $m(s,t)\neq 2$, and having components
corresponding to direct product factors of $W$. The Coxeter graph
is not used in this paper.) The presentation diagram of the
special subgroup of $W$ generated by a subset $S'\subseteq S$ is
the induced subgraph of $\Gamma(W,S)$ with vertex set $S'$ (and
in this sense, special subgroups could as well be termed {\it
visual} subgroups since we can see the presentation diagram of
such a subgroup in $\Gamma(W,S)$).

Suppose $\Gamma(W,S)=\Gamma_1\cup\Gamma_2$ is a union of induced
subgraphs and let $\Gamma_0=\Gamma_1\cap\Gamma_2$ (so vertices and
edges of $\Gamma(W,S)$ are in $\Gamma_1$ or $\Gamma_2$ or both,
and $\Gamma_0$ is the induced subgraph consisting of the vertices
and edges in both).  Equivalently, suppose $\Gamma_0$ is an
induced subgraph with $\Gamma(W,S)-\Gamma_0$ having at least two
components, $\Gamma_1$ is $\Gamma_0$ together with some of these
components and $\Gamma_2$ is $\Gamma_0$ together with the other
components. We say in this case that $\Gamma_0$ {\it separates}
$\Gamma(W,S)$ (separates it into at least two components).  Then
it is evident from the Coxeter presentation that $W$ is an
amalgamated product of special subgroups corresponding to
$\Gamma_1$ and $\Gamma_2$ over the special subgroup corresponding
to $\Gamma_0$. Amalgamated product decompositions with special
factors and special amalgamated subgroup are easily seen in the
presentation  diagram and we call such an amalgamated product a
{\it visual} splitting of $W$. Other amalgamated product
decompositions may also be possible, and we want to understand
such splittings in terms of visual splittings.

More generally, we are interested in when a Coxeter group $W$ can
be realized as the fundamental group of a graph of groups (as
explained in the next section). We show that the graph must
actually be a tree and so this generalizes amalgamated products.
We say that $\Psi$ is a {\it visual graph of groups decomposition}
of $W$ (for a given $S$), if each vertex and edge group of $\Psi$
is a special subgroup of $W$, the injections of each edge group
into its endpoint vertex groups are given simply by inclusion, and
the fundamental group of $\Psi$ is isomorphic to $W$ by the
homomorphism induced by the inclusion map of vertex groups into
$W$.  A sequence of compatible visual splittings of $W$ will
result in such a decomposition. Our main result shows that an
arbitrary graph of group decomposition of a Coxeter group can be
refined (in a certain sense) to a visual graph of groups
decomposition.

\begin{theorem}[Main Theorem]\label{maintheorem}
Suppose $(W,S)$ is a Coxeter system and $W$ is a subgroup of the
fundamental group of a graph of groups $\Lambda$. Then $W$ has a
visual graph of groups decomposition $\Psi$ where each vertex
group of $\Psi$ is a subgroup of a conjugate of a vertex group of
$\Lambda$, and each edge group of $\Psi$ is a subgroup of a
conjugate of an edge group of $\Lambda$.  Moreover, $\Psi$ can be
taken so that each special subgroup of $W$ that is a subgroup of
a conjugate of a vertex group of $\Lambda$ is a subgroup of a
vertex group of $\Psi$.
\end{theorem}

For $(W,S)$ and $\Lambda$ as in theorem
\ref{maintheorem} and $\Psi$ satisfying the full conclusion of theorem \ref{maintheorem} (including the moreover clause),  we say $\Psi$ is a $(W,S)$-{\it visual
decomposition from} $\Lambda$.  Suppose $G$ is a group decomposed
as $A\ast _CB$. If $H$ is a subgroup of $B$, then $\langle A\cup
H\rangle$ decomposes as $A\ast _C \langle C\cup H\rangle$
(consider the action of $\langle A\cup H\rangle$ on the
Bass-Serre tree for $A\ast _CB$). Furthermore, the decomposition
$(A\ast _C \langle C\cup H\rangle )_{\langle C\cup H\rangle} B$
reduces to $A\ast _C B$ so $G$ decomposes as $\langle A\cup
H\rangle \ast _{\langle C\cup H\rangle} B$. In this way some
visual decompositions of Coxeter groups can be ``artificially"
altered to decompositions that have no visual vertex or edge
groups. In section \ref{maintheoremproof}, we exhibit a Coxeter
system $(W,S)$ and a reduced graph of groups decomposition
$\Lambda _N$, for $(W,S)$ with $N$ vertices for any positive
integer $N$. Even so, the following theorem defines limits on how
far an arbitrary graph of groups decomposition for a finitely
generated Coxeter system can stray from a visual decomposition
for that system.

\begin{theorem}\label{artificial}
Suppose $(W,S)$ is a finitely generated Coxeter system, $\Lambda$
is a graph of groups decomposition of $W$ and $\Psi$ is a reduced
graph of groups decomposition of $W$ such that each vertex group of $\Psi$ is a subgroup of a conjugate of a vertex group of $\Lambda$. Then for each vertex
$V$ of $\Lambda$, the vertex group $\Lambda (V)$, has a graph of
groups decomposition $\Phi_V$ such that each vertex group of
$\Phi _V$ is either

(1) conjugate to a vertex group of $\Psi$ or

(2) a subgroup of $v\Lambda (E)v^{-1}$ for some $v\in \Lambda
(V)$ and $E$ some edge of $\Lambda$ adjacent to $V$.
\end{theorem}

When $\Psi$ is visual, vertex groups of the first type in theorem \ref{artificial} are
visual. Those of the second type seem somewhat artificial.

It is easy to recognize whether a finitely generated Coxeter
group is 2-ended or infinite ended, by refining Stallings'
theorem (from \cite{Stallings}) to a visual splitting theorem.  A
Dunwoody decomposition of a finitely generated Coxeter group is a
graph of groups with finite or 1-ended vertex groups and finite
edge groups. These refine to a visual Dunwoody decomposition, and
we get a simple argument for why finitely generated Coxeter
groups are accessible with respect to splittings over finite
groups (as in \cite{Dunwoody}). In separate papers
\cite{MTACCESS} and \cite{MTJSJ}, we developed the fundamental
ideas of this paper to prove a strong accessibility result for
Coxeter groups with respect to splittings over ``minimal"
splitting subgroups, and a JSJ result for splittings of Coxeter
groups over virtually abelian groups. At the bottom level, the
strong accessibility result is a visual version of Dunwoody's
result in the Coxeter setting.

J. P. Serre gives an account of FA groups in \cite{Serre}.  In
particular, an FA group has no nontrivial splittings. We apply
the main theorem to show that the maximal FA subgroups of a
finitely generated Coxeter group $W$ are those conjugate to a
special subgroup with presentation diagram a maximal complete
subgraph of $\Gamma(W,S)$.

Our results identify  certain properties of a Coxeter group
recognizable directly from a particular presentation of the
group, properties apparent in the presentation  diagram. There
can be different Coxeter systems $(W,S)$ and $(W,S')$ where $S'$
is not simply a conjugate of $S$ (and may not even correspond
under any automorphism of $W$), and so conjugates of special
subgroups with respect to $S$ need not correspond to conjugates
of special subgroups with respect to $S'$.  However, these
results do imply that certain special subgroups are, up to
conjugation, special for any Coxeter system.  This investigation
thus has significant application to rigidity questions and the
determination of when Coxeter groups are isomorphic.

As a final applications of these ideas we give a visual
classification of finitely generated virtually free Coxeter
groups in section \ref{applicationssection}.

The main theorem generalizes to settings of groups other than
Coxeter groups, e.g., for graph products of finite groups.  Our
arguments work for essentially any graph presentation with finite
groups on vertices and edges, satisfying a certain
``developability'' hypothesis. Even more generally, the
hypothesis of finite groups in a graph presentation is only used
to ensure that these groups lie in a conjugate of a vertex group
of the given graph of groups decomposition, and so very similar
arguments apply to graph products with infinite vertex and edge
groups provided these groups are assumed to lie in conjugates of
the vertex groups of the given graph of groups decomposition.

\section{Graphs of Groups}\label{gogsection}

Our main tool in this investigation is the connection between
group actions on trees and fundamental groups of graphs of
groups.  As reference the reader is referred to
\cite{DicksDunwoody} and \cite{Serre}.  We review some of the
pertinent definitions and results.

A graph of groups $\Lambda$ consists of a set $V(\Lambda)$ of
vertices, a set $E(\Lambda)$ of edges, and maps
$\iota,\tau:E(\Lambda)\to V(\Lambda)$ giving the initial and
terminal vertices of each edge in a connected graph, together with
vertex groups $\Lambda(V)$ for $V\in V(\Lambda)$, edge groups
$\Lambda(E)$ for $E\in E(\Lambda)$, with $\Lambda(E)\subset
\Lambda(\iota(E))$ and an injective group homomorphism
$t_E:\Lambda(E)\to\Lambda(\tau(E))$, called the edge map of $E$
and denoted by $t_E:g\mapsto g^{t_E}$. The fundamental group
$\pi_1(\Lambda)$ of a graph of groups $\Lambda$ is the group with
presentation having generators the disjoint union of $\Lambda(V)$
for $V\in V(\Lambda)$, together with a symbol $t_E$ for each edge
$E\in E(\Lambda)$, and having as defining relations the relations
for each $\Lambda(V)$, the relations $gt_E=t_Eg^{t_E}$ for $E\in
E(\Lambda)$ and $g\in\Lambda(\iota(E))$, and relations $t_E=1$
for $E$ in a given spanning tree of $\Lambda$ (the result, up to
isomorphism, is independent of the spanning tree taken).

An amalgamated product $A*_CB$ is realized as the fundamental
group of the graph of groups with 2 vertices having vertex groups
$A$ and $B$, a single edge between, with edge group (the image in
$A$ of) $C$, and edge map $t_E$ determined by the injection of the
edge group into $B$. Similarly, an HNN-extension of $A$ by
$t_E:C\to A$ is realized as the fundamental group of a graph of
groups with a single vertex and single edge. In general, the
fundamental group of a graph of groups can be understood as taking
amalgamated products of the vertex groups along the edge groups in
the spanning tree, followed by HNN-extensions over the remaining
edge groups.  The $t_E$ for edges not in the spanning tree
correspond to a stable letter of an HNN extension. For edges in
the spanning tree the relations amount to identifying $\Lambda(E)$
in $\Lambda(\iota(E))$ with its image $\Lambda(E)^{t_E}$ in
$\Lambda(\tau(E))$ as in an amalgamated product.  Each vertex and
edge group of a graph of groups $\Lambda$ injects into the
fundamental group of $\Lambda$ (Britton's lemma) and we usually
identify these with the corresponding subgroups of the fundamental
group of $\Lambda$.

A graph of groups on a graph which is not simply connected will
have $\mathbb Z$ as a homomorphic image (take a $t_E$ not in the
spanning tree to 1 and all other generators to 0).  Since the
generators of a Coxeter group are of order 2, a homomorphism into
$\mathbb Z$ must have trivial image.  Thus if a Coxeter group is
the fundamental group of a graph of groups the graph must be a
tree and the group arises as successive amalgamated products.  In
working with such trees of groups we will often simply assume
suitable identifications have been made and the edge maps are
simply inclusion maps.

Given a graph of groups $\Lambda$, the Bass-Serre tree for
$\Lambda$ is defined with vertices the disjoint union over $V\in
V(\Lambda)$ of the different cosets $g\Lambda(V)$ of $\Lambda(V)$
in $\pi_1(\Lambda)$, and edges the disjoint union over $E\in
E(\Lambda)$ of the different cosets $g\Lambda(E)$ of $\Lambda(E)$
in $\pi_1(\Lambda)$, taken with
$\iota(g\Lambda(E))=g\Lambda(\iota(E))$ and
$\tau(g\Lambda(E))=g\Lambda(\tau(E))$.  The Bass-Serre tree is in
fact a tree and the fundamental group of $\Lambda$ acts on this
tree by taking for $h\in \pi_1(\Lambda)$,
$h(g\Lambda(V))=(hg)\Lambda(V)$ and
$h(g\Lambda(E))=(hg)\Lambda(E)$.

If a group $G$ acts on a tree $T$ (as a directed graph, the action
preserving the orientation of edges), then a transversal for this
action consists of a vertex and edge from each orbit of the action
of $G$ on vertices and edges.  There must always exist a
transversal having a spanning subtree such that each other edge of
the transversal originates in the subtree. (E.g., if $\Lambda$ has
a single vertex and single edge, giving $G=\pi_1(\Lambda)$ as an
HNN extension, then $G$ acts on the Bass-Serre tree $T$ for
$\Lambda$ with transversal consisting of a single vertex and an
edge originating at that vertex, that does not include the
terminal vertex, and has spanning subtree just the single
vertex.) From such a transversal, a graph of groups is defined by
taking as graph the quotient of $T$ under the action of $G$, and
taking vertex and edge groups to be the stabilizers of the
corresponding vertices and edges in the transversal.  For edges
in the spanning subtree of the transversal the edge maps are
given by inclusion. An edge of the transversal not in the
spanning subtree connects a vertex in the transversal to a
translate by an element $g\in G$ of a vertex in the transversal
and the edge map is given by conjugation by $g$.  Then the
fundamental group of this graph of groups is naturally isomorphic
to $G$ by the homomorphism extending the inclusion map of vertex
groups into $G$.  Choosing a different transversal in $T$ gives
rise to a different graph of groups decomposition of $G$ with
vertex groups each conjugate to vertex groups of the first graph
of groups, i.e., having isomorphic vertex groups but
corresponding with different subgroups of $G$.

Suppose $\Lambda$ is a graph of groups and $X$ is a subset of the
edges in $\Lambda$.  Contracting the graph of $\Lambda$ along the
edges in $X$ gives rise to a graph whose vertices are the
equivalence classes of vertices modulo the equivalence relation
defined by identifying endpoints of edges in $X$, taking edges
$E\in E(\Lambda)-X$ with endpoints the equivalence classes of
endpoints of $E$.  Define a graph of groups $\Lambda'$ on this
graph, having the same fundamental group as $\Lambda$, by taking
as vertex groups $\Lambda'(V)$ the fundamental group of the graph
of groups with vertices of $\Lambda$ identified to $V$ and edges
in $X$ between these vertices, with edge groups $\Lambda'(E)$
corresponding to the remaining edge groups in $\Lambda$ with
corresponding edge maps.  We call $\Lambda'$ the graph of groups
resulting from collapsing the edges $X$ in $\Lambda$.  If
$\Lambda$ is a tree of groups, then collapsing edges in $\Lambda$
results in another tree of groups.

A graph of groups will be called {\it reduced} if no edge group
is equal to its originating vertex group nor has image equal to
its terminating vertex group, EXCEPT for an edge which is a loop
at a vertex. If a graph of groups is not reduced, then we may
collapse a vertex across an edge, where the edge group is (or has
image) the same as the endpoint vertex group, giving a smaller
graph of groups with vertex and edge groups among the original
vertex and edge groups and having the same fundamental group.
Repeated reductions of this sort in a finite graph of groups must
eventually end with a reduced graph of groups all of whose vertex
and edge groups were present in the original graph of groups and
having the same fundamental group.

The following is well-known. We include a proof as this result is
frequently referenced throughout the paper.

\begin{lemma} \label{R1}
Suppose $\Lambda$ is a reduced graph of groups decomposition of a
group $G$, the underlying graph for $\Lambda$ is a tree, $V$ and
$U$ are vertices of $\Lambda$, and $g\Lambda(V)g^{-1}\subset
\Lambda(U)$ for some $g\in G$, then $V=U$ and $g\in\Lambda(V)$.
\end{lemma}
\begin{proof}
If $U\ne V$ or if $U=V$ and $g\not\in \Lambda (V)$, then in the
Bass-Serre tree for $\Lambda$, $\Lambda(V)$ stabilizes the
distinct vertices $\Lambda(V)$ and $g^{-1}\Lambda(U)$.  But then
$\Lambda(V)$ also stabilizes a geodesic path between these
vertices and hence stabilizes the first edge $h\Lambda(E)$ for an
edge $E$ of $\Lambda$ at $V$. This would mean $\Lambda(V)$ is
equal to $\Lambda(E)$ (or its image in $\Lambda(V)$), and we
could collapse $V$ across $E$, contradicting that $\Lambda$ was
reduced.
\end{proof}

In the other direction, given a graph of groups $\Lambda$ and a
graph of groups decomposition of a vertex group $\Lambda(V)$ as
$\pi_1(\Phi)$, we would like to see when $\Lambda$ results from
collapsing out $\Phi$ in a larger graph of groups $\Lambda'$. Say
that $\Phi$ is a {\it compatible decomposition} of $\Lambda(V)$ if
each edge group of $\Lambda$ incident at $V$ is a subgroup of a
$\Lambda(V)$ conjugate of a vertex group of $\Phi$. In general
then, given $\Lambda$ and $\Phi$ a compatible decomposition of
$\Lambda(V)$, construct $\Lambda'$ by replacing $V$ in $\Lambda$
by $\Phi$ and attaching edges of $\Lambda$ incident at $V$ instead
to vertices in $\Phi$ with edge groups and edge maps appropriately
adjusted so each such edge group has image in a vertex group of
$\Phi$ (for edges ending in $\Phi$ in $\Lambda'$ or $V$ in
$\Lambda$) or is an isomorphic subgroup of a vertex group of
$\Phi$ (for edges originating in $\Phi$ in $\Lambda'$ or $V$ in
$\Lambda$). This adjustment of edge maps means that the image in
$\pi_1(\Lambda')$ of a vertex group of $\Lambda'$ not in $\Phi$ is
conjugate to the image of the corresponding vertex group of
$\Lambda$ in $\pi_1(\Lambda)=\pi_1(\Lambda')$.

Alternatively, the Bass-Serre tree $T'$ for $\Lambda'$ can be
constructed from the Bass-Serre tree for $\Lambda$ with each coset
of $\Lambda(V)$ replaced by the Bass-Serre tree for $\Phi$ with a
$\pi_1(\Lambda)$ group action.  Contracting the edges in the orbit
of $\Phi$ in $T'$ gives $T$, and a transversal of $T'$ with a
spanning subtree contracts to a transversal of $T$ with a spanning
subtree. Starting with a transversal with spanning subtree of $T$
and a transversal with spanning subtree of the Bass-Serre tree for
$\Phi$ we may need to translate the parts of the $T$ transversal
in different components of $T-\lbrace V\rbrace$ so that they
attach to the transversal of the Bass-Serre tree for $\Phi$ in
$T'$.  The result is still a transversal of $T$ with the same
quotient graph as $\Lambda$ but with vertex groups identified with
conjugates of the vertex groups of $\pi_1(\Lambda)$.

Suppose $(W,S)$ is a Coxeter system.  A visual graph of groups
decomposition $\Psi$ of $W$ has special vertex and edge groups,
edge maps given by inclusion, and is such that the inclusion of
vertex groups in $W$ extends to an isomorphism of $\pi_1(\Psi)$
with $W$. As noted above, $\Psi$ is a tree.  In terms of
presentation diagrams, $\Gamma(W,S)$ must be the union of the
subdiagrams corresponding to the vertex groups of $\Psi$, with
the edge groups corresponding to the intersections of adjacent
vertex subdiagrams.  To understand in visual terms exactly when a
graph of special subgroups has fundamental group isomorphic to
$W$ we have the following essential lemma.

\begin{lemma}\label{visgog}
Suppose $(W,S)$ is a Coxeter system.  A graph of groups $\Psi$
with graph a tree, where each vertex group and edge group is a
special subgroup and each edge map is given by inclusion, is a
visual graph of groups decomposition of $W$ iff each edge in the
presentation diagram of $W$ is an edge in the presentation diagram
of a vertex group and, for each generator $s\in S$, the set of
vertices and edges with groups containing $s$ is a nonempty
subtree in $\Psi$.
\end{lemma}

\begin{proof}
Suppose $\Psi$ is a visual graph of groups decomposition of $W$ so
that the homomorphism $\psi$ extending the inclusion map on vertex
groups is an isomorphism of $\pi_1(\Psi)$ and $W$.  Since each
vertex group is generated by the elements of $S$ it contains, the
image of $\psi$ is generated by the set of all $s\in S$ that
belong to some vertex group.  Since a proper subset of $S$
generates a proper subgroup of $W$, each element of $S$ must be an
element of a vertex group of $\Psi$. Let $T$ be the Bass-Serre
tree for $\Psi$. Then since the edge maps are inclusions, the
subgraph of $T$ consisting of the identity cosets of vertex and
edge groups of $\Psi$ is a transversal of $T$, a subtree of $T$,
and the stabilizers of the vertices and edges of the transversal
are simply the corresponding vertex and edge groups of $\Psi$.  If
$s\in S$ belongs to two vertex groups, then $\psi^{-1}(s)$
stabilizes the corresponding vertices of the transversal, hence
stabilizes all the edges and vertices in a geodesic path in $T$
between these vertices, and so $s$ belongs to all the vertex and
edge groups in the path in $\Psi$ between these vertices.  Hence
the vertices and edges with groups containing $s$ form a subtree
in $\Psi$. If there is an edge in $\Gamma(W,S)$ between $s,t\in
S$, then the subgroup $U$ of $W$ generated by $s$ and $t$ is
finite. But then $\psi^{-1}(U)$ is a finite subgroup of
$\pi_1(\Psi)$ acting on $T$ and so must stabilize a vertex $V$ of
$T$.  Let $V'$ be the vertex of the transversal closest to $V$ in
$T$.  Then $V'$ lies between $V$ and a vertex of the transversal
stabilized by $s$ and so is also stabilized by $s$, and similarly
for $t$.  Thus $V'$ corresponds to a vertex of $\Psi$ having
vertex group containing both $s$ and $t$ and this special subgroup
has an edge between $s$ and $t$ in its presentation diagram as
well.

Conversely, suppose $\Psi$ is a visual graph of groups with each
edge of the presentation  diagram of $W$ in the presentation
diagram of some vertex group, and such that, for each generator
$s\in S$, the vertices and edges of $\Psi$ with groups containing
$s$ form a subtree in $\Psi$. All of the occurrences of an $s\in
S$ in different vertex groups of $\Psi$ are identified by the
relators of $\pi_1(\Psi)$ along edges in the subtree of vertex and
edge groups containing $s$.  Take $\psi:\pi_1(\Psi)\to W$
extending the inclusion map of the generators $S$ in $\pi_1(\Psi)$
into $W$. We get that $\psi$ is an isomorphism by checking that
each of the defining relators of $W$ is already a relator of
$\pi_1(\Psi)$. But each generator of $W$ belongs to some vertex
group of $\Psi$ where it has order 2, and the other relators of
$W$ correspond to edges of the presentation  diagram and these
relators also already hold in some vertex group of $\Psi$.
\end{proof}

If $\Psi$ is a visual graph of groups decomposition for the
Coxeter system $(W,S)$, it is convenient to label a vertex of $\Psi$ by the subset of $S$ that generates the
corresponding vertex group. So if $Q\subset S$ is a
vertex label of $\Psi$, then $\Psi (Q)=\langle Q\rangle$.
If $\Psi$ is reduced, its vertex labels are distinct and
we identify vertices with their labels. Even if $\Psi$ is reduced, two distinct edges may have the same edge group, so we do not extend this labeling to edges.

The following two technical results are useful tools.

\begin{corollary} \label{C1}
Suppose $(W,S)$ is a finitely generated Coxeter system, $\Psi$ is a
visual graph of groups decomposition of $W$, and $E\subset S$ is such that $\langle E\rangle=\Psi (E')$ for $E'$ an edge of $\Psi$. If $\{x,y\}\subset S-E$, and $x\in X$ and
$y\in Y$ for $X$ and $Y$ labels of vertices of $\Psi$ on opposite
sides of $E'$, then as a subset of $\Gamma(W,S)$, $E$ separates
$x$ and $y$ in $\Gamma$.
\end{corollary}

\begin{proof} Otherwise, among all such $x,y$ that fail the
conclusion, let $\alpha$ be a shortest path in $\Gamma$ from $x$
to $y$ avoiding $E$. Note that by lemma \ref{visgog}, $x\ne y$. If
$x\equiv x_0,x_1,\ldots x_n\equiv y$ are the consecutive vertices
of $\alpha$, then by lemma \ref{visgog} $\{x_0,x_1\}\subset V$
for some vertex label $V$ of $\Psi$. As $x_0\not\in E$, lemma
\ref{visgog} implies $V$ and $X$ label vertices on the same side
of $E'$ in $\Psi$. But then $x_1$ and $y$ satisfy the hypothesis
of the corollary contradicting the minimality of $\alpha$.
\end{proof}

\begin{corollary}\label{C2}
Suppose $(W,S)$ is a finitely generated Coxeter system, and
$\Psi$ is a visual graph of groups decomposition of $W$. If $C$ is
a complete subset of the presentation diagram $\Gamma (W,S)$, then
there is a vertex labeled $V(\subset S)$ of $\Psi$ such that
$C\subset V$.
\end{corollary}
\begin{proof} We may assume $\Psi$ is reduced. If $C$ is not a
subset of a vertex of $\Psi$, then there is a vertex $U\subset S$
of $\Psi$ containing a maximal number of elements of $C$ and
$V\subset S$ a vertex of $\Psi$ closest to $U$ such that $V$
contains an element $c$ of $C-U$. If $E'$ is the last edge of the
$\Psi$-geodesic from $U$ to $V$ and $\Psi(E')=\langle E\rangle$ for $E\subset S$, then $(U\cap C)\not \subset E$
and $c\not \in E$. By corollary \ref{C1}, $E$ separates $(U\cap
C)-E$ from $c$ in $\Gamma$. This is impossible as $C$ is complete.
\end{proof}

\section{Proof of Main Results}\label{maintheoremproof}

We begin this section with a proof of theorem \ref{maintheorem} and conclude with a proof of theorem \ref{artificial}. Several examples are presented to introduce the reader to visual decompositions.

\medskip

\noindent {\bf Proof of Theorem \ref{maintheorem}:} 
Suppose $(W,S)$ is a Coxeter system and suppose $\Lambda$ is a
graph of groups with $W$ a subgroup of $\pi_1(\Lambda)$.  We may
identify the vertex and edge groups of $\Lambda$ with subgroups
of $\pi_1(\Lambda)$.  Then $\pi_1(\Lambda)$ acts on the Bass-Serre
tree $T$ such that the stabilizers of vertices of the tree are
conjugates of the vertex groups of $\Lambda$ and stabilizers of
edges of the tree are conjugates of the edge groups of $\Lambda$.

We build a visual graph of groups decomposition $\Psi$ of $W$ with
this tree as its graph. For each vertex $V$ (resp. edge $E$) of
$\Psi$ take vertex group $\Psi(V)$ (resp. edge group $\Psi(E)$) to
be the subgroup of $W$ generated by the $s\in S$ stabilizing $V$
(resp. $E$).  The edge groups inject into the vertex groups by
inclusion maps. Clearly, each vertex group (resp. edge group) of
$\Psi$ is a special subgroup of $W$ and a subgroup of a conjugate
of a vertex group (resp. edge group) of $\Lambda$.  Each generator
$s\in S$ is of order 2 and so stabilizes some vertex of $T$.  If
two vertices of $T$ are stabilized by $s$, then $s$ stabilizes the
geodesic path between the vertices, thus the subgraph of $T$ of
vertices and edges stabilized by $s$ is a subtree.  Suppose there
is an edge in the presentation diagram for $(W,S)$ between $s$ and
$t$ in $S$. The subgroup generated by $s$ and $t$ has finite
order, so stabilizes some vertex of $T$, with $s$ and $t$
belonging to that vertex group. By Lemma \ref{visgog} then,
$\pi_1(\Psi)$ is isomorphic to $W$ and $\Psi$ is a visual graph of
groups decomposition of $W$.

Finally, if $U$ is a special subgroup and a subgroup of a
conjugate of one of the vertex groups of $\Lambda$, then $U$
stabilizes a vertex in $T$ and so $U$ is a subgroup of the vertex
group of that vertex.  $\square$

\smallskip
Generally, $\Psi$ will be an infinite tree.  If $W$ is finitely
generated, by taking the subtree spanning a finite set of vertices
whose vertex groups between them include all of the elements of
$S$, we get a finite visual graph of groups decomposition of $W$.
Reducing this graph of groups results in a reduced visual graph of
groups decomposition for $W$.

\begin{remark}\label{R2} If $W$, $\Lambda$ and $\Psi$ are as in theorem
\ref{maintheorem}, then each edge $E$ of $\Lambda$ defines a
splitting of $W$ as $A\ast _{\Lambda (E)}B$ where $A$ is
generated by the vertex groups of the vertices on one side of $E$
in $\Lambda$ and $B$ is generated by those on the other side of
$E$ in $\Lambda$. There is a quotient map from the Bass-Serre tree
for $\Lambda$ to that of $A\ast _{\Lambda (E)}B$ which respects
the action of $W$. If $A\not =\Lambda (E)\not =B$ then $W$ is not
a subgroup of a conjugate of $A$ or $B$ and the proof above shows
that there is an edge of $\Psi$ such that the edge divides $\Psi$
into two components, neither alone having all of the generators
of $W$ in its vertex groups, and such that the corresponding edge
group lies in a conjugate of $\Lambda (E)$ (and so also after
applying possible reductions to $\Psi$).
\end{remark}

\begin{example}\label{simex}
As an example of the main theorem consider the Coxeter group $W$
with presentation
$$P=\langle s_{1},s_{2},s_{3},s_{4},s_{5}:\,s_{i}^{2},\,
(s_{1}s_{2})^{2}=(s_{2}s_{3})^{2}=(s_{3}s_{4})^{2}=(s_{4}s_{5})^{2}=1\rangle$$
Note that $(s_{3}s_{5}s_{3})^{2}=(s_{4}s_{3}s_{5}s_{3})^{2}=1$ so
the map $\phi:W\to W$ defined on generators by $\phi(s_{i})=s_{i}$
for $i<5$, and $\phi(s_{5})=s_{3}s_{5}s_{3}$, extends to a
homomorphism, is surjective, is its own inverse, and so is an
automorphism of $W$.  Consider subgroups $A=\langle
s_{1},s_{2},s_{4},s_{3}s_{5}s_{3}\rangle$, $B=\langle
s_{2},s_{3},s_{4}\rangle$, and $C=\langle s_{2},s_{4}\rangle$.
Then $W=A*_{C}B$ is the image under $\phi$ of a visual amalgamated
product decomposition of $W$, but $A$ is not a visual subgroup of
$W$.  We refine this decomposition to a visual graph of groups
decomposition following the proof of theorem \ref{maintheorem}.

The Bass-Serre tree on which $W$ acts has vertices corresponding
to the different cosets of $A$ and $B$ and edges corresponding to
the different cosets of $C$ with $hC$ connecting $hA$ and $hB$ for
each $h\in W$.  The vertex $A$ (respectively $B$) of this tree has
one edge for each coset $aC$ (resp. $bC$) of $C$ in $A$ (resp.
$B$). An element $g\in W$ acts on this tree by mapping $hA$ to
$ghA$ and similarly for cosets of $B$ and $C$. The stabilizer of
$gA$ is the conjugate $gAg^{-1}$. The graph of groups $\Psi$ is
defined on this tree by taking the vertex and edge groups
generated by the $s_{i}$ which stabilize that vertex or edge. Thus
$s_{1}$, $s_{2}$, and $s_{4}$ stabilize the $A$ vertex in the
tree, $s_{2}$ and $s_{4}$ stabilize both the $C$ edge and the
$s_{3}C$ edge, $s_{2}$, $s_{3}$, and $s_{4}$ stabilize the
$B=s_{3}B$ vertex, $s_{2}$, $s_{4}$, and $s_{5}$ stabilize the
$s_{3}A$ vertex (etc.). Since the vertex groups of these three
vertices include all of the generators of $W$, the vertex and edge
groups for the rest of the tree collapse back to this three vertex
subgraph with the same fundamental group as the graph of groups.
The visual decomposition reduces to a three factor amalgamated
product
$$W=\langle s_{1},s_{2},s_{4}\rangle*_{C}B*_{C}\langle
s_{2},s_{4},s_{5}\rangle$$ where the first and third factors are
generated by those generators of $W$ stabilizing $A$ and $s_{3}A$,
respectively, are subgroups of conjugates of $A$, and are special
subgroups of $W$.

This simple example illustrates that an amalgamated product
decomposition of a Coxeter group need not refine to a visual
amalgamated product decomposition of only two factors.  Instead
the visual graph of groups decomposition we have produced is
sometimes necessary.
\end{example}

\noindent {\bf Proof of Theorem \ref{artificial}:} Let $\Phi _V$ be the
decomposition of $\Lambda(V)$ given by selecting a fundamental
transversal for its action on $T_{\Psi}$, the Bass-Serre tree for
$\Psi$.  That is, take the graph of $\Phi_V$ to be a subtree of
$T_{\Psi}$ containing a single vertex and edge from each orbit of
the action of $\Lambda(V)$ on $T_{\Psi}$ and take the group of a
vertex or edge to be the subgroup of $\Lambda(V)$ stabilizing
that vertex or edge.  If $U$ is a vertex of $\Phi_V$, then $U\in
T_{\Psi}$ corresponds to a coset $g\Psi(X)$ of a vertex group
$\Psi(X)$ of $\Psi$. The stabilizer in $W$ of $U$ is simply the
conjugate $g\Psi(X)g^{-1}$, and so $\Phi_V(U)=\Lambda(V)\cap
(g\Psi(X)g^{-1})$. The group $\Psi(X)$ is a subgroup
of a conjugate $h\Lambda(Y)h^{-1}$ of a vertex group $\Lambda(Y)$
of $\Lambda$. Hence $\Phi_V(U)$ is a subgroup of $\Lambda(V)$ and
of a conjugate $gh\Lambda(Y)h^{-1}g^{-1}$ of $\Lambda(Y)$, and so
stabilizes the vertices of $T_{\Lambda}$ corresponding to
$\Lambda(V)$ and $gh\Lambda(Y)$.  If $V=Y$ and $gh\in\Lambda(V)$,
then these are the same vertex of $T_{\Lambda}$, and
$g\Psi(X)g^{-1}\subset gh\Lambda(Y)h^{-1}g^{-1}=\Lambda(V)$ so
$\Phi_V(U)=\Lambda(V)\cap (g\Psi(X)g^{-1})=g\Psi(X)g^{-1}$ is
equal to a conjugate of a vertex group of $\Psi$.  Otherwise,
$\Lambda(V)$ and $gh\Lambda(Y)$ are different vertices of
$T_{\Lambda}$. The group $\Phi_V(U)$ stabilizes both vertices and
so a geodesic path between the two. In particular, $\Phi_V(U)$
stabilizes the first edge $k\Lambda(E)$ in this path, where $E$ is
an edge of $\Lambda$ incident with $V$ and $k\in\Lambda(V)$, and
thus $\Phi_V(U)\subset k\Lambda(E)k^{-1}$. Reducing $\Phi_V$
gives the desired graph of groups decomposition. $\square$

\medskip

The next corollary and example examine technical aspects of the
decompositions of theorem \ref{artificial}.

\begin{corollary}\label{artificialcor}
Suppose $(W,S)$ is a finitely generated Coxeter system, $\Lambda$
is a graph of groups decomposition of $W$, $\Psi$ is a reduced
$(W,S)$-visual decomposition from $\Lambda$ and $\Phi_V$ is a
reduced graph of groups decomposition for the vertex group
$\Lambda (V)$ as given in the proof of theorem \ref{artificial}.
If $X$ is a vertex of $\Psi$, $V$ is a vertex of $\Lambda$, and
$g\Psi(X)g^{-1}$ is a subgroup of $\Lambda(V)$ for some $g\in W$,
then $vg\Psi(X)g^{-1}v^{-1}$ is a vertex group of $\Phi_V$ for
some $v\in \Lambda (V)$.
\end{corollary}

\begin{proof} Let $\Phi_V'$ be the decomposition of $\Lambda (V)$
given by selecting a fundamental transversal for its action on
$T_{\Psi}$, the Bass-Serre tree for $\Psi$, and let $\Phi_V$ be
reduced from $\Phi_V'$. By the definition of transversal, there
is $v\in \Lambda (V)$ such that the coset $vg\Psi (X)$
corresponds to a vertex $B$ of the transversal. Then
$\Phi_V'(B)=\Lambda (V)\cap
vg\Psi(X)g^{-1}v^{-1}=vg\Psi(X)g^{-1}v^{-1}$. It remains to show
that this vertex group survives reduction. Otherwise, there is a
vertex $Q$ of $\Phi_V'$ such that $vg\Psi(X)g^{-1}v^{-1}$ is a
proper subgroup of $\Phi_V'(Q)\equiv \Lambda(V)\cap
h\Psi(Z)h^{-1}$ (where $Q\in T_{\Psi}$ corresponds to the coset
$h\Psi(Z)$). But then $\Psi(X)$ is a proper subgroup of
$g^{-1}v^{-1}h\Psi(Z)h^{-1}vg$ which is impossible by lemma
\ref{R1}.
\end{proof}

The next example shows that the previous corollary cannot be
extended to show that a vertex group of $\Psi$ is conjugate to a
vertex group of $\Phi_V$ for a unique vertex $V$ of $\Lambda$. In
fact, all vertex groups of a particular $\Phi_V$ may be equal to
vertex groups of other $\Phi_Q$ for $Q$ a vertex of $\Lambda$.

\begin{example}\label{artificialex}
Consider the group $W$ with Coxeter presentation $\langle w,x,y,z:
w^2=x^2=y^2=z^2=1\rangle$ (and free product decomposition
$\langle w\rangle\ast \langle x\rangle \ast\langle y\rangle\ast
\langle z\rangle$). There is an automorphism of $W$ that fixes $w$
and $y$, and maps $x$ to $wyxyw$ and $z$ to $y(wyxyw)z(wyxyw)y$.
Observe that $W$ has the decomposition $\langle w,x\rangle
\ast_{\langle x\rangle} \langle x,z\rangle\ast_{\langle z\rangle}
\langle y,z\rangle$. Hence the automorphism induces a
decomposition $\Lambda$, with vertices $V_1$ and $V_2$ and $V_3$,
where $V_1$ and $V_2$ are connected by an edge $E_1$, and $V_2$
and $V_3$ are connected by an edge $E_2$. The decomposition
$\Lambda\equiv \Lambda(V_1)\ast _{\Lambda(E_1)}\Lambda(V_2)\ast
_{\Lambda (E_2)}\Lambda (V_3)$ of $W$ is such that

$\Lambda (V_1)=\langle w,wyxyw\rangle$

$\Lambda (E_1)=\langle wyxyw\rangle$

$\Lambda(V_2)= \langle wyxyw, y(wyxyw)z(wyxyw)y\rangle$

$\Lambda(E_2)=\langle y(wyxyw)z(wyxyw)y\rangle$ and

$\Lambda(V_3)= \langle y, y(wyxyw)z(wyxyw)y \rangle$.

A visual decomposition $\Psi$, from $\Lambda$ is $\langle w\rangle
\ast \langle x\rangle \ast \langle y\rangle \ast \langle
z\rangle$. The decomposition $\Phi_{V_1}$ of theorem
\ref{artificial} is $\langle w\rangle \ast \langle
wyxyw\rangle$,  $\Phi_{V_2}$ is $\langle wyxyw\rangle\ast \langle
y(wyxyw)z(wyxyw)y\rangle$, and $\Phi_{V_3}$ is $\langle y\rangle
\ast \langle y(wyxyw)z(wyxyw)y\rangle$. Hence the vertex groups
$\langle x\rangle$ and $\langle z\rangle$ of $\Psi$ are conjugate
to vertex groups of $\Phi_V$ for more than one vertex $V$ of
$\Lambda$. Furthermore, both vertex groups of $\Phi_{V_2}$ are
equal to vertex groups of other $\Phi_{V_i}$.
\end{example}

\begin{example}\label{HNN}
Consider the Coxeter presentation $\langle a,b,c :
a^2=b^2=c^2=1\rangle$. This group splits as $\Lambda\equiv
\langle a,bc\rangle\ast _{\langle bc\rangle}\langle b,c\rangle$.
The visual decomposition from this splitting is $\Psi\equiv
\langle a\rangle \ast \langle b,c\rangle$. If $\Lambda (V)$ is
the vertex group $\langle a,bc\rangle$, then its graph of groups
decomposition $\Phi_V$, induced by its action on $T_{\Psi}$, the
Bass-Serre tree for $\Psi$, is $\mathbb Z_2\ast \mathbb Z$. This
group also has a decomposition with underlying graph not tree,
but it is not induced by its action on $T_{\Psi}$. We wonder if
every decomposition of a finitely generated Coxeter group and
corresponding visual decomposition can only induce $\Phi_V$ have
underlying graph a tree.
\end{example}

In their accessibility paper \cite{BF}, Bestvina and Feighn limit
the number of vertex groups that can occur in a reduced graph of
groups decomposition with ``small" edge groups, for a given
finitely presented group. Strong accessibility results for Coxeter
groups over ``minimal" splittings are the focus of
\cite{MTACCESS}. While theorem \ref{artificial} limits how far a
vertex group of an arbitrary decomposition of a Coxeter group
can stray from visual, the following example shows that one
cannot expect an accessibility result for Coxeter groups over
arbitrary splittings.

\begin{example}
For a virtually free (non-virtually cyclic) group, the following
technique is standard for creating non-trivial reduced graph of
groups decompositions with $n$ vertex groups for any $n\in
\{1,2,\ldots \}$. Let $W$ be the Coxeter group with presentation
$\langle s_1,s_2,s_3,s_4,s_5,s_6: s_i^2=1\rangle$. Then $W$ has
free product decomposition $A\ast B$ where $A=\langle
s_1,s_2,s_3\rangle$ and $B=\langle s_4,s_5\rangle$. The element
$b\equiv s_4s_5$ generates an infinite cyclic subgroup of finite
index in $B$ and $A$ contains a non-abelian free subgroup of
finite index. Let $\{a_1,a_2,\ldots \}$ be a free generating set
for an infinite rank free subgroup of $A$. Then $W$ has the
decomposition
$$\Lambda _1=\langle A\cup\{b\}\rangle \ast _{\langle
a_1,b\rangle} \langle B\cup \{a_1\}\rangle$$ The group $\langle
A\cup \{b\}\rangle$ has the following non-trivial reduced graph
of groups decomposition, which is compatible with $\Lambda_1$ :
$$\langle A \cup \{b^2\}\rangle \ast _{\langle
a_1,a_2, b^2\rangle} \langle a_1,a_2,b\rangle$$ Hence $W$ has the
reduced decomposition: $$\Lambda_2=\langle A \cup \{b^2\}\rangle
\ast _{\langle a_1,a_2,b^2\rangle} \langle a_1,a_2,b\rangle \ast
_{\langle a_1,b\rangle} \langle B\cup \{a_1\}\rangle$$ Defining
$E_n=\langle a_1,\ldots ,a_{n},b^{(2^{n-2})}\rangle$, $C_n=\langle
a_1,\ldots ,a_{n},b^{(2^{n-1})}\rangle $ and $A_n=\langle A\cup
\{b^{(2^{n-1})}\}\rangle$, we have a reduced decomposition for $W$
with $n+1$ vertex groups:
$$\Lambda_n= A_n \ast _{C_n} E_n \ast _{C_{n-1}} E_{n-1}
\ast _{C_{n-2}}\ast \ldots \ast E_2\ast _{C_1}\langle B\cup
\{a_1\}\rangle$$ For each $n$, $A\ast B$ is a visual
decomposition from $\Lambda _n$.
\end{example}

\section{Technical Results \label{Tech}}

As mentioned in the introduction, our basic reference for Coxeter
groups is \cite{Bourbaki}. The technical results of this section
are well-known and can be derived from two fundamental facts: the
first concerning special subgroups of Coxeter groups and the
second called the ``deletion condition".

\begin{proposition}\label{flat}
Suppose $(W,S)$ is a Coxeter system and $P=\langle S:
(st)^{m(s,t)}$ for $m(s,t)<\infty\rangle$ (where $m:S^2\to
\{1,2,\ldots ,\infty\}$) is a Coxeter presentation for $W$. If
$A\subset S$, then $(\langle A\rangle, A)$ is a Coxeter system
with Coxeter presentation $\langle A:(st)^{m'(s,t)}$ for
$m'(s,t)<\infty\rangle$ (where $m'=m\vert _{A^2}$).
\end{proposition}

Given a group $G$ and a generating set $S$, an {\it $S$-geodesic}
for $g\in G$ is a shortest word in $S\cup S^{-1}$ such that the
product of the letters of this word is $g$. The number of letters
in an $S$-geodesic for $g$ is the {\it $S$-length of $g$}.  We
include an elementary proof of the deletion condition pointed out
to us by A. Y. Ol'shanskii, and based on the theory of van Kampen
diagrams (see chapter 5 of Lyndon and Schupp's book \cite{LS} for
a basic introduction to van Kampen diagrams).

\begin{proposition}\label{dc} {\bf (The Deletion Condition)} Suppose
$(W,S)$ is a Coxeter system and $w=a_1\cdots a_n$ for $a_i\in S$.
If $a_1\cdots a_n$ is not geodesic then there are indices $i<j$
in $\{1,2,\ldots ,n\}$ such that $w=a_1\cdots a_{i-1}a_{i+1}\cdots
a_{j-1}a_{j+1}\cdots a_n$. I.e. $a_i$ and $a_j$ can be deleted.
\end{proposition}

\begin{proof}
Let $b_1\cdots b_m$ be an $S$-geodesic for $w$, so that $m<n$. Let
$\cal D$ be a van Kampen diagram for the word $a_1\cdots
a_nb_m\cdots b_1$. The relation 2-cell $r_1$ containing $a_1$ is
of even length and so there is an edge $e_1$ of $r_1$ ``opposite"
$a_1$. The two subpaths of the boundary of $r_1$ separated by
$a_1$ and $e_1$ have the exact same labeling. Let $r_2$ be the
relation 2-cell of $\cal D$ that shares $e_1$ with $r_1$, and let
$e_2$ be the edge of $r_2$ opposite $e_1$. Again the labeling of
the two subpaths of the boundary of $r_2$ separated by $e_1$ and
$e_2$ are the same. Continue until $e_k$ is on the boundary of
$\cal D$. The ``strip" determined $\cup _{i=1}^kr_i$ is such that
the two subpaths of the boundary of this strip separated by $a_1$
and $e_k$ have the same labeling. This creates a unique paring of
edges on the boundary of $\cal D$. Since $m<n$ there are indices
$i<j$ in $\{1,2,\ldots ,n\}$ such that $a_i$ is paired with
$a_j$. If the two subpaths of the boundary of the strip for
$a_i$, separated by $a_i$ and $a_j$ are labeled $\beta$, then (as
the product of edge path labels of loops in $\cal D$ are trivial
in $W$) $a_{i+1}\cdots a_{j-1}=\beta =a_i\cdots a_j$. I.e. $a_i$
and $a_j$ delete.
\end{proof}

\begin{lemma}\label{lett}
Suppose $(W,S)$ is a Coxeter system, $A\subset S$, and $a\in
\langle A\rangle$. If $a=a_1\cdots a_n=b_1\cdots b_n$ are
$S$-geodesics then $\{a_1,\ldots ,a_n\}=\{b_1,\ldots ,b_n\}$.
\end{lemma}

\begin{proof}
Assume $n$ is minimal among all counterexamples to the lemma. By
proposition \ref{flat}, $n\ne 1$. Note that $a_1\cdots
a_nb_n=b_1\cdots b_{n-1}$. By the deletion condition, $b_n$
deletes with some $a_i$ in the first expression. Inductively
$\{b_1,\ldots ,b_{n-1}\}\subset \{a_1,\ldots ,a_n\}$. Similarly,
$b_1a_1\cdots a_n=b_2,\cdots b_n$ so that $\{b_2,\ldots
,b_n\}\subset \{a_1,\ldots ,a_n\}$.
\end{proof}

If $G$ is a group with generating set $S$, then the {\it Cayley
graph} ${\cal K}(G,S)$ is a labeled graph with vertex set $G$ and
a directed edge (labeled $s$) from the vertex $g$ to the vertex
$gs$ for each $s\in S$. Given a vertex $x$ of ${\cal K}(G,S)$,
there is a bijective correspondence between edge paths at $x$ and
words in $S\cup S^{-1}$, where traversing an edge labeled $s$
opposite its orientation is read as $s^{-1}$. Hence $S$-geodesics
and ${\cal K}(G,S)$-geodesics at a given vertex are the same.

The next result is a straightforward application of the deletion
condition.

\begin{lemma}\label{concat}
Suppose $(W,S)$ is a Coxeter system, $\{x,y\}\subset
W$ and $A\subset S$. Then in the Cayley graph ${\cal K}(W,S)$
there is a unique closest point $z$ of the coset $y\langle
A\rangle$ to $x$. Furthermore, if $\alpha$ is the $S$-geodesic
from $x$ to $z$, and $\beta$ is any $A$-geodesic at $z$, then
$\alpha \beta$ is an $S$-geodesic.
\end{lemma}

The next lemma follows by a result of Kilmoyer (see Theorem 2.7.4
of \cite{Carson}), but we include a direct proof for completeness.

\begin{lemma}\label{intersectspecial}
Suppose $(W,S)$ is a Coxeter system, $I,J\subseteq S$, and $d$ is
a minimal length double coset representative in $\langle I\rangle
d \langle J\rangle$.  Then $\langle I\rangle\cap d\langle
J\rangle d^{-1}=\langle K\rangle$ where $K=I\cap(dJd^{-1})$.
Hence $g\langle I\rangle g^{-1}\cap h\langle J\rangle
h^{-1}=f\langle K\rangle f^{-1}$ for $K\subset I$.
\end{lemma}

\begin{proof}
With $K=I\cap(dJd^{-1})$, clearly, $\langle K\rangle \subset
\langle I\rangle\cap d\langle J\rangle d^{-1}$.  Suppose $a$ is a
shortest element in $\langle I\rangle\cap d\langle J\rangle
d^{-1}$ but not in $\langle K\rangle$.  Write $a=a_1a_2\ldots a_m$
geodesically with $a_i\in I$ and $d^{-1}ad=b_1b_2\ldots b_n$
geodesically for $b_i\in J$. Write $d=d_1d_2\ldots d_k$
geodesically for $d_i\in S$. Then since $d$ is a minimal length
double coset representative, $ad=a_1a_2\ldots a_md_1d_2\ldots
d_k=d_1d_2\ldots d_kb_1b_2\ldots b_n$ are each geodesic.  Hence
$m=n$. Clearly $a\neq 1$, and $a\notin S$ else
$a=a_1=db_1d^{-1}\in K$ by definition.  Instead, $m>1$.  Now
$b_1b_2\ldots b_md^{-1}$ and $d^{-1}a_m$ are geodesic but
$d^{-1}a_1\ldots a_{m-1}=b_1b_2\ldots b_md^{-1}a_m$ so this last
is not geodesic and, by the deletion condition, $a_m$ deletes with
some $b_i$ to give $b_1b_2\ldots b_md^{-1}a_m=b_1\ldots
b_{i-1}b_{i+1}\ldots b_md^{-1}$ (if $1<i<m$ and similarly if $i=1$
or $i=m$).  But then $a_1a_2\ldots a_{m-1}=db_1\ldots
b_{i-1}b_{i+1}\ldots b_md^{-1}\in \langle I\rangle\cap d\langle
J\rangle d^{-1}$ and, by the minimality of $a$ and lemma
\ref{lett}, we have $\{a_1,a_2,\ldots ,a_{m-1}\}\subset K$.
Likewise $\{a_m,\ldots,a_2\}\subset \langle I\rangle\cap d\langle
J\rangle d^{-1}$ so $a_m\in K$.  But then $a\in\langle K\rangle$
contradicting the choice of $a$. Instead, every $a\in\langle
I\rangle\cap d\langle J\rangle d^{-1}$ is in $\langle K\rangle$.

Given conjugates $g\langle I\rangle g^{-1}$ and $h\langle
J\rangle h^{-1}$ of special subgroups, take $d$ of minimal length
in $\langle I\rangle g^{-1}h\langle J\rangle$ so
$g^{-1}h=adb^{-1}$ for $a\in\langle I\rangle$ and $b\in\langle
J\rangle$, $\langle I\rangle g^{-1}h\langle J\rangle=\langle
I\rangle d\langle J\rangle$, and $\langle I\rangle\cap d\langle
J\rangle d^{-1}=\langle K\rangle$ as above.  Then

\begin{eqnarray*}
g\langle I\rangle g^{-1}\cap h\langle J\rangle h^{-1}&=&ga\langle
I\rangle a^{-1}g^{-1}\cap hb\langle J\rangle b^{-1}h^{-1}\\
&=&ga(\langle I\rangle\cap d\langle J\rangle d^{-1})a^{-1}g^{-1}\\
&=&ga\langle K\rangle a^{-1}g^{-1}
\end{eqnarray*}
\end{proof}

We will use the following corollary in section \ref{FAsection} in
our analysis of FA subgroups.

\begin{corollary}\label{specialconj}
Suppose $(W,S)$ is a finitely generated Coxeter system,
$I,J\subseteq S$, the induced subgraph on $I$ is a maximal
complete subgraph of $\Gamma(W,S)$, and $I\subseteq w\langle
J\rangle w^{-1}$ for some $w\in W$. Then $I\subseteq J$ and
$w\in\langle J\rangle$.
\end{corollary}

\begin{proof}
Take $d$ a minimal length double coset representative in $\langle
I\rangle w\langle J\rangle$ so $w=adb^{-1}$ for $a\in\langle
I\rangle$ and $b\in\langle J\rangle$. Then since
\begin{eqnarray*}
\langle I\rangle\cap d\langle J\rangle d^{-1}&=& a^{-1}(a\langle
I\rangle a^{-1}\cap adb^{-1}\langle J\rangle bd^{-1}a^{-1})a\\
&=&a^{-1}(\langle I\rangle\cap w\langle J\rangle w^{-1})a\\
&=&\langle I\rangle
\end{eqnarray*}
we have $I\subseteq dJd^{-1}$ so $d^{-1}Id\subseteq J$.  If $d\neq
1$ and $d=d_1d_2\ldots d_k$ geodesically, then by the minimality
of $d$, $d_1\notin I$. There is some $s\in I$  unrelated to $d_1$,
since $I$ induces a maximal complete subgraph of the Coxeter
diagram. But if $a_1\ldots a_k$ is any geodesic and $s$ is
unrelated to $a_k$ then $a_1\ldots a_ks$ is geodesic (apply lemma
\ref{concat} with $x=1$, $y=a_1\cdots a_{k-1}$, $A=\{a_k,s\}$ and
$\beta =\beta _1a_ks$, where $\beta _1$ is the $\{a_k,s\}$
geodesic from $z$ to $y$). Hence $d^{-1}sd$ cannot have length 1
and cannot be in $J$. Instead then, $d=1$, $I\subseteq J$ and
$w=adb^{-1}\in\langle J\rangle$.
\end{proof}

The following corollaries find application in the analysis of
rigidity of Coxeter groups.

\begin{corollary}\label{Icontained}
If $(W,S)$ is a finitely generated Coxeter system, $I\subseteq S$,
and $w^{-1}\langle I\rangle w\subset \langle I\rangle$, for some
$w\in W$, then $w^{-1}\langle I\rangle w=\langle I\rangle$.
\end{corollary}
\begin{proof}
Let $d$ be a minimal length double coset representative in
$\langle I\rangle w\langle I\rangle$, so $w=adb^{-1}$ for
$a,b\in\langle I\rangle$.  Then
\begin{eqnarray*}
\langle I\rangle&=&\langle I\rangle\cap w\langle I\rangle w^{-1}\\
&=&a^{-1}(\langle I\rangle\cap d\langle I\rangle d^{-1})a
\end{eqnarray*}
Hence $\langle I\rangle=\langle I\rangle\cap d\langle I\rangle
d^{-1}=\langle I\cap dId^{-1}\rangle$, by lemma
\ref{intersectspecial}. By propositions \ref{flat} and \ref{dc}, a
proper subset of $I$ cannot generate the same Coxeter group as
$I$, so $I=dId^{-1}$. Hence $d^{-1}Id=I$ and $w^{-1}\langle
I\rangle w=bd^{-1}\langle I\rangle db^{-1}=b\langle I\rangle
b^{-1}=\langle I\rangle$.
\end{proof}

\begin{corollary}
Suppose $(W,S)$ and $(W,S')$ are finitely generated Coxeter
systems for the same Coxeter group $W$.  Suppose $I\subset S$ is
such that the induced subgraph of $I$ separates $\Gamma(W,S)$ (its
complement has at least two components).  Then there are sets
$J\subset S$ and $J'\subset S'$ such that the induced subgraphs
of $J$ and $J'$ separate $\Gamma(W,S)$ and $\Gamma(W,S')$,
respectively, $\langle J\rangle$ and $\langle J'\rangle$ are
conjugate, and $\langle J\rangle$ is conjugate to a subgroup of
$\langle I\rangle$.
\end{corollary}

\begin{proof}
As noted in the introduction, since $I$ separates $\Gamma(W,S)$,
$W=A*_{\langle I\rangle}B$, for $A$ the ($S$) special subgroup
generated by $I$ and the generators in some components of the
complement of $I$, and $B$ the ($S$) special subgroup generated by
$I$ and the generators in the other components of the complement
of $I$ in $\Gamma(W,S)$, $A\neq\langle I\rangle\neq B$.  Taking
this amalgamated product decomposition of $W$ to be $\Lambda$, we
consider a corresponding visual decomposition $\Psi$ with respect
to the alternate generating set $S'$.  As noted in remark
\ref{R2}, there must be an edge of $\Psi$, with not all of the
generators $S'$ appearing on one side, having edge group $\langle
I'_1\rangle$, for an $I'_1\subset S'$, which is a subgroup of a
conjugate of $\langle I\rangle$. That is, there is an
$I'_1\subset S'$ which separates $\Gamma(W,S')$ and for which
$\langle I'_1\rangle$ is a subset of a conjugate of $\langle
I\rangle$.

But then $I'_1$ gives a proper splitting of $W$ and by the same
reasoning there is an $I_2\subset S$ which separates
$\Gamma(W,S)$ and which is a subset of a conjugate of $\langle
I'_1\rangle$. Continuing in this fashion, we can find
$I'_{2k+1}\subset S'$ separating $\Gamma(W,S')$, generating a
subgroup of a conjugate of $\langle I_{2k}\rangle$, and an
$I_{2k+2}\subset S$ separating $\Gamma(W,S)$ and generating a
subgroup of a conjugate of $\langle I'_{2k+1}\rangle$.

Since there are only finitely many subsets of $S$ or $S'$, for
some $k_1<k_2<k_3$, $I'_{k_1}=I'_{k_3}$. Take
$J'=I'_{k_1}=I'_{k_3}$ and $J=I_{k_2}$.  Since $\langle
I_j\rangle$ and $\langle I'_{j+1}\rangle$ are subgroups of
conjugates of $\langle I_i\rangle$ and $\langle I'_{i+1}\rangle$
for $i<j$, take $a_1,a_2\in W$ such that
$$
\langle I'_{k_3}\rangle \subset a_2\langle I_{k_2}\rangle
a_2^{-1} \subset a_2a_1\langle I'_{k_1}\rangle a_1^{-1}a_2^{-1}
$$
Then $a_2^{-1}a_1^{-1}\langle J'\rangle a_1a_2 \subset \langle
J'\rangle$, so by corollary \ref{Icontained} these are equal, and
we get
$$
\langle J'\rangle \subset a_2\langle J\rangle a_2^{-1} \subset
a_2a_1\langle J'\rangle a_1^{-1}a_2^{-1}=\langle J'\rangle
$$
so $\langle J'\rangle= a_2\langle J\rangle a_2^{-1}$ and $J$ and
$J'$ generate subgroups of conjugates of $\langle I\rangle$ that
separate $\Gamma(W,S)$ and $\Gamma(W,S')$ respectively.
\end{proof}

\section{Ends}\label{oddsandends}

Stalling's theorem \cite{Stallings} states that if a finitely
generated group has more than one end then it splits nontrivially
as an amalgamated product or HNN-extension over a finite group.
The following result is then an easy consequence of our main
theorem.  It can also be obtained from work of M. Davis
\cite{Davis}.

\begin{corollary}\label{corends}
For any finitely generated Coxeter group $W$ with presentation
diagram $\Gamma$, the following are equivalent
\begin{enumerate}
\item $W$ has more than one end
\item $W$ decomposes as a nontrivial amalgamated product $A*_{C}B$
where $C$ is finite and $A$, $B$, and $C$ are special subgroups
\item $\Gamma$ contains a complete separating subgraph, the vertices
of which generate a finite subgroup of $W$.
\end{enumerate}
\end{corollary}

\begin{proof}
If $W$ is not 1-ended or finite, then by Stallings' splitting
theorem, $W=A*_CB$ with $C$ finite (and not as an HNN-extension
since an HNN-extension maps onto $\mathbb Z$ but a homomorphism of
a Coxeter group into $\mathbb Z$ must take generators to the
identity).  Then by theorem \ref{maintheorem}, $W$ has a reduced
visual graph of groups decomposition in which each edge group is a
subgroup of a conjugate of $C$ and so is finite.  This
decomposition cannot be trivial since each vertex group is a
subgroup of a conjugate of $A$ or $B$ neither of which is $W$ in
the given nontrivial splitting, and so no vertex group can be $W$.
Hence there is at least one edge after reducing $\Psi$ and
collapsing the other edges gives $W$ as a visual amalgamated
product of special subgroups over a finite special subgroup.  The
remaining implications are easy.
\end{proof}

By Stallings' theorem and theorem \ref{maintheorem}, a 2-ended
Coxeter group splits as a visual amalgamated product over a
finite group which is of index two in each factor, the following
result then characterizes 2-ended Coxeter groups.

\begin{corollary}
A Coxeter group with system $(W,S)$ and diagram $\Gamma$ is
2-ended iff $\Gamma$ contains a separating subdiagram
$\Gamma_{0}$ which is the presentation diagram of a finite group,
and $\Gamma-\Gamma_{0}$ consists of two vertices each of which is
connected to each vertex of $\Gamma_{0}$ by edges labeled 2 (but
not connected to each other). Equivalently, $W= \langle x,y\rangle
\times \langle H\rangle $ where $\{x,y\}\cup H= S$, $x$ and $y$
are unrelated and $\langle H\rangle$ is finite.
\end{corollary}

Thus the number of ends of a Coxeter group can be easily
determined from an analysis of separating subdiagrams of a Coxeter
diagram, and checking which subdiagrams correspond to finite
subgroups. A Coxeter group whose presentation diagram is complete
is either finite or 1-ended. The finite Coxeter groups have been
enumerated \cite{Bourbaki}.

A Dunwoody decomposition of a finitely presented group is a graph
of groups decomposition of the group with finite edge groups and
1-ended and finite vertex groups.  In \cite{Dunwoody}, Dunwoody
shows that any finitely presented group has such a decomposition.

\begin{corollary}\label{DDproposition}
Suppose $(W,S)$ is a finitely generated Coxeter system and $W$ is
the fundamental group of a graph of groups $\Lambda$ where each
edge group is finite.  Then $W$ has a visual decomposition $\Psi$
where each vertex group is 1-ended or finite and a subgroup of a
conjugate of a vertex group of $\Lambda$, and where each edge
group is finite.
\end{corollary}

\begin{proof}
Take $\Psi$ the reduced visual graph of groups from $\Lambda$ as
given by theorem \ref{maintheorem}.  Then each edge group of
$\Psi$ is a subgroup of a conjugate of an edge group of $\Lambda$
and so is finite.  Suppose some vertex $V(\subset S)$ of $\Psi$
is such that $\langle V\rangle$ is not 1-ended or finite.  Then
$\langle V\rangle$ visually splits nontrivially over a finite
subgroup. If $E(\subset S)$ is an edge of $\Psi$ with endpoint
$V$, then $E\subset V$ and $\langle E\rangle$ is finite. In
particular, $E$ induces a complete subgraph in $\Gamma(W,S)$. By
corollary \ref{C2}, $E$ must be contained in a vertex group of
the visual decomposition of $\langle V\rangle$, i.e., the
splitting of $\langle V\rangle$ is visually compatible with
$\Psi$. Replacing $\langle V\rangle$ in $\Psi$ by this splitting
gives a visual graph of groups decomposition with finite edge
groups. Since a special vertex group is replaced by special
vertex groups with fewer generators, repeating this process
eventually must end with a visual graph of groups decomposition
having finite edge groups and finite or 1-ended vertex groups.
\end{proof}

One might wonder how ``visual'' a Dunwoody decomposition must be.

\begin{theorem}\label{closedunwoody}
Suppose $(W,S)$ is a finitely generated Coxeter system. If
$\Lambda$ is a reduced Dunwoody graph of groups decomposition of
$W$ and $\Psi$ is a reduced visual decomposition for $(W,S)$ such
that each edge group of $\Psi$ is finite and each vertex group of
$\Psi$ is a subgroup of a conjugate of a vertex group of $\Lambda$
(in particular if $\Psi$ is a reduced visual graph of groups
decomposition from $\Lambda$), then
\begin{enumerate}
\item  $\Psi$ is a Dunwoody decomposition

\item  There is a (unique) bijection $\alpha$ of the vertices
of $\Lambda$ to the vertices of $\Psi$ such that for each vertex
$V$ of $\Lambda$, $\Lambda(V)$ is conjugate to $\Psi(\alpha (V))$

\item  Each edge group of $\Lambda$ is conjugate to a special
subgroup for $(W,S)$.
\end{enumerate}
\end{theorem}

\begin{proof}
If $\Phi$ is a graph of groups decomposition of $W$ with finite
edge groups, then any finite or 1-ended subgroup of $W$ is a
subgroup of a conjugate of a vertex group of $\Phi$ (otherwise,
the action of this group on the Bass-Serre tree for $\Phi$ would
induce a non-trivial splitting over a finite group). Hence each
vertex group of $\Lambda$ is a subgroup of a conjugate of a
vertex group of $\Psi$.  If a vertex group $A=\Lambda(V)$ of
$\Lambda$ is a subgroup of a conjugate of $\Psi(V')$ for $V'$ a
vertex of $\Psi$, then since $\Psi(V')$ is a subgroup of a
conjugate of a vertex group of $\Lambda$, $A$ is a subgroup of a
conjugate of a vertex group of a vertex $V''$ of $\Lambda$.  As
noted in lemma \ref{R1}, in a reduced graph of groups, the vertex
group $A$ at $V$ is a subgroup of a conjugate of a vertex group
at $V''$ only if $V=V''$ and the conjugate is by an element of
$A$. But then $A$ is conjugate to $\Psi(V')$. Again, since no
vertex group of $\Psi$ is contained in a conjugate of another,
$V'$ is uniquely determined, and we set $\alpha(V)=V'$. Since
each vertex group $\Psi(V')$ is contained in a conjugate of some
$\Lambda(V)$ which is in turn conjugate to $\Psi(\alpha(V))$ we
must have $V'=\alpha(V)$ and each $V'$ is in the image of
$\alpha$.  In particular, each vertex group of $\Psi$ is 1-ended
or finite and so $\Psi$ is a Dunwoody decomposition of $W$.

Since $\Lambda$ is a tree, we can take each edge group of
$\Lambda$ as contained in its endpoint vertex groups taken as
subgroups of $W$. Hence each edge group is simply the intersection
of its adjacent vertex groups (up to conjugation). Since vertex
groups of $\Lambda$ correspond to conjugates of vertex groups in
$\Psi$, their intersection is conjugate to a special subgroup by
Lemma \ref{intersectspecial}.
\end{proof}

In Example \ref{simex} we have a visual Dunwoody decomposition
$$\langle s_1,s_2\rangle*_{\langle s_2\rangle}
\langle s_2,s_3\rangle*_{\langle s_3\rangle} \langle
s_3,s_4\rangle*_{\langle s_4\rangle} \langle s_4,s_5\rangle
$$
which is carried by the automorphism $\phi$ to the Dunwoody
decomposition where the last factor is replaced by $\langle
s_4,s_3s_5s_3\rangle$.  Thus in this theorem we cannot expect a
single element to conjugate all factors of a Dunwoody
decomposition to the corresponding factors of the corresponding
visual decomposition.  The connection between an arbitrary
Dunwoody decomposition and a visual Dunwoody decomposition is
however clearly quite close.

It is worthwhile to see how this analysis of Coxeter groups leads
to an understanding of why finitely generated Coxeter groups are
accessible. While this argument only re-proves a special case of
Dunwoody's accessibility theorem, it is the base case of the main
theorems of our papers \cite{MTACCESS} and \cite{MTJSJ}, where we
prove a strong accessibility result for Coxeter groups and
splittings over ``minimal" splitting subgroups, and a JSJ result
for Coxeter groups and splittings over virtually abelian
subgroups.

\begin{lemma}
Suppose $(W,S)$ is a finitely generated Coxeter system and
$\Lambda$ is a graph of groups decomposition of $W$ with finite
edge groups.  Suppose a vertex group of $\Lambda$ splits
nontrivially as $A*_CB$ over a finite $C$. Then there is a special
subgroup or a subgroup of a finite special subgroup of $W$
contained in a conjugate of $B$ which is not also contained in a
conjugate of $A$ (and then also with $A$ and $B$ reversed).
\end{lemma}

\begin{proof}
Let $\Lambda'$ be the graph of groups resulting from replacing the
vertex whose group splits by the graph corresponding to $A*_CB$.
Let $\Psi$ be the corresponding visual graph of groups
decomposition of $W$ and let $T$ be the Bass-Serre tree for
$\Psi$. The intersection of any conjugates of $A$ and $B$, or the
intersection of distinct conjugates of $B$, is contained in a
conjugate of an edge group and so is finite. If an infinite vertex
group of $\Psi$ lies in a conjugate of $B$, it cannot also lie in
a conjugate of $A$, so suppose no infinite vertex group of $\Psi$
lies in a conjugate of $B$.  From the action of $B$ on $T$ we get
a reduced graph of groups decomposition of $B$ with vertex and
edge groups contained in conjugates of vertex and edge groups of
$\Psi$, in particular all of the edge groups are finite.  If any
vertex group $B_1$ of this decomposition is infinite, it is
contained in a conjugate of an infinite vertex group of $\Psi$
which is in turn contained in a conjugate of a vertex group of
$\Lambda'$ other than $B$.  But $B_1$ would then be an infinite
subgroup of $B$ and contained in a conjugate of another vertex
group of $\Lambda'$ which is impossible. Instead, the vertex
groups of the decomposition of $B$ are finite and conjugate to
subgroups of finite special subgroups. Replace $B$ in $\Lambda'$
by this graph of groups decomposition and collapse the edge $C$ if
it equals one of the vertex groups of the decomposition of $B$ to
get a new graph of groups decomposition $\Lambda''$ where $A$ is
adjacent to a vertex group $B_1$ of the decomposition of $B$ by an
edge labeled $C_1$, a proper subgroup of $B_1$ (either with
$C_1=C$, if no collapse happens, or with $C_1$ an edge group of
the decomposition of $B$). Now $B_1$ is finite, contained in $B$,
but $B_1$ cannot also be contained in a conjugate of $A$, since
otherwise $B_1$ would stabilize a path from a coset of $A$ to
$B_1$ in the Bass-Serre tree for $\Lambda''$ and hence stabilize a
coset of $C_1$, i.e., would be contained in a conjugate of $C_1$
which has fewer elements than $B_1$.
\end{proof}

As mentioned earlier, our next theorem follows from Dunwoody's
accessibility theorem, and the proof of this theorem leads to a
general approach to more complex accessibility and JSJ results.

\begin{theorem}\label{Dunacc}
Finitely generated Coxeter groups are accessible.
\end{theorem}

\begin{proof}
Suppose $(W,S)$ is a finitely generated Coxeter system. There are
only finitely many special subgroups of $W$ and finitely many
subgroups of finite special subgroups.  For $G$ a subgroup of $W$
let $n(G)$ be the number of special subgroups or subgroups of
finite special subgroups which are contained in any conjugate of
$G$ (which includes the trivial group), so $1\leq n(G)\leq n(W)$.
For $\Lambda$ a finite graph of groups decomposition of $W$ let
$c(\Lambda)=(c_{n(W)},\ldots,c_2,c_1)$ where $c_i$ is the count of
vertex groups $G$ of $\Lambda$ with $n(G)=i$.  Let $<$ order these
$n(W)$-tuples lexicographically, a well ordering of $n(W)$-tuples
of nonnegative integers.  If $\Lambda$ reduces to $\Lambda'$ then
clearly no $c_i$ increases and some $c_i$ must decrease.  If a
vertex group $G$ of $\Lambda$ splits as $A*_CB$ to produce a new
$\Lambda'$, then every subgroup of a conjugate of $A$ or $B$ is a
subgroup of a conjugate of $G$, but, by the last lemma, some
special subgroup or subgroup of a finite special subgroup is
contained in a conjugate of $B$, and so of $G$, but not in a
conjugate of $A$. Hence $n(A)<n(G)$, and similarly $n(B)<n(G)$,
and so $c(\Lambda')<c(\Lambda)$ since $c_{n(G)}$ decreases by 1 in
going from $\Lambda$ to $\Lambda'$ and the only other components
that change are $c_{n(A)}$ and $c_{n(B)}$ which are later in the
tuples. Since $<$ is a well ordering, there can be no infinite
sequence of graph of group decompositions of $W$ resulting from
successive reductions or splittings over finite subgroups, i.e.,
$W$ is accessible.
\end{proof}

\section{Maximal FA subgroups of Coxeter Groups}\label{FAsection}

Let $(W,S)$ be a finitely generated Coxeter system. A set of
vertices of a complete subgraph of $\Gamma (W,S)$ is called a {\it
simplex} of $(W,S)$ and the subgroup of $W$ generated by a simplex
of $(W,S)$ is called a {\it simplex group} of $(W,S)$. If $V$ is a
vertex of a graph $\Gamma$, define $lk(V)$, {\it the link of $ V$}, to be the set of all
vertices $V'$ of $\Gamma$ such that $V'$ is connected to $V$ by an
edge. Define $st(V)$, {\it the star of $V$}, to be $lk(V)\cup \{V\}$.

A group $G$ is called $FA$ if for every tree on which $G$ acts,
the set of fixed points of $G$ in the tree is non-empty. In $\S6$
of \cite{Serre}, we find basic results on FA groups. In
particular we find:

\begin{proposition}\label{S1}
If $G$ is denumerable then  $G$ is FA iff $G$ is finitely
generated, no quotient of $G$ is isomorphic to the infinite
cyclic group and $G$ does not split non-trivially as an
amalgamated product.
\end{proposition}

\begin{proposition}\label{S2}
If an FA group $G$ is a subgroup of $A*_{C}B$ then $G$ is a
subgroup of a conjugate of $A$ or $B$.
\end{proposition}

\begin{proposition}\label{Every}
Every simplex subgroup of a Coxeter system is FA.
\end{proposition}

We have then the following result.

\begin{lemma}\label{L3.1} Suppose $(W,S)$ is a Coxeter system and $G$ is an
FA-subgroup of W, then $G$ is a subgroup of a conjugate of a
simplex subgroup of $(W,S)$.
\end{lemma}

\begin{proof} If the presentation  diagram $\Gamma(W,S)$
is not complete, choose vertices $s$ and $t$ that are not
related. Then $lk(s)$ separates $s$ and $t$ in $\Gamma(W,S)$.
Hence

$$W\cong \langle st(s)\rangle \ast _{\langle lk(s)\rangle
}\langle S-\{s\}\rangle$$

By proposition \ref{S2}, $G$ is a subgroup of a conjugate of
$\langle st(s)\rangle $ or $\langle S-\{s\}\rangle$. Continue
splitting until the conclusion is realized.
\end{proof}

The following theorem is the main result of this section. Our
proof is based in group actions on trees. A more combinatorial
approach works equally well.

\begin{theorem}\label{maxfatheorem}
The maximal FA subgroups of a Coxeter group $W$ are precisely the
conjugates of the special subgroups whose diagrams are the maximal
complete subdiagrams of a (and equivalently any) presentation
diagram for $W$
\end{theorem}

\begin{proof}
We have that simplex subgroups are FA and any FA subgroup is
contained in a conjugate of a simplex subgroup.  If $A$ is a
maximal simplex subgroup contained in $wBw^{-1}$ for $B$ another
simplex subgroup, then by corollary \ref{specialconj},  $A=B$ and
$w\in B$.
\end{proof}

By a result of Tits \cite{Bourbaki}, up to conjugacy there are
only finitely many elements of order 2 in a finitely generated
Coxeter group. Hence we have:

\begin{proposition} Suppose $(W,S)$ is a finitely generated Coxeter
system. Then the subgroup of automorphisms $a$ of $W$ such that
$a(s)$ is a conjugate of $s$ for all $s\in S$, is of finite index
in $Aut(W)$. In particular, there exists an integer $n(W,S)$ such
that for any $a\in Aut(W)$, $a^{n}(s)$ is conjugate to $s$ for
all $s\in S$.
\end{proposition}

Next we prove:

\begin{theorem} Suppose $(W,S)$ is a finitely generated Coxeter
system with maximal simplices $\sigma _{1},\ldots ,\sigma _{m}$.
If $C$ is the subgroup of all $c\in Aut(W)$ such that, for
$i\in \{1,\ldots ,m\}$ there exists a $w_{i,c}\in W$ so that, $c(x)=w_{i,c}xw_{i,c}^{-1}$
when $x\in \langle \sigma _{i}\rangle $ (i.e., $c$ restricted to
$\langle\sigma_i\rangle$ is conjugation by a $w_{i,c}\in W$
depending only on $c$ and $\sigma _{i}$), then $C$ has finite
index in $Aut(W)$.
\end{theorem}

\begin{proof}
If $a\in Aut(W)$ and $\sigma _{1},\ldots ,\sigma _{m}$ are the
maximal simplicies of $(W,S)$, then $a(\langle \sigma _{i}\rangle
)=w_{i,a}\langle \sigma _{\alpha (i)}\rangle w_{i,a}^{-1}$ for some
$\alpha (i)\in \lbrace 1,\ldots ,m\rbrace$, by theorem
\ref{maxfatheorem}. Observe that $\sigma _{\alpha (j)}\not =\sigma
_{\alpha (k)}$ for $k\not =j$, since otherwise
$$
a^{-1}(w_{k,a}^{-1})\langle \sigma _{k}\rangle a^{-1}(w_{k,a})\equiv
a^{-1}(\langle \sigma _{\alpha (k)}\rangle)=a^{-1}(\langle\sigma
_{ \alpha (j)}\rangle)\equiv a^{-1}(w_{j,a}^{-1})\langle \sigma
_{j}\rangle a^{-1}(w_{j,a})
$$
implies $\sigma_k=\sigma_j$ by corollary \ref{specialconj}. Hence
$\alpha $ is a permutation of $(1,\ldots ,m)$. If $b\in Aut(W)$
and $b(\langle \sigma _{j}\rangle )=w_{j,b}\langle \sigma _{\beta
(j)}\rangle w_{j,b}^{-1}$ then 
$$ba(\langle \sigma _{i}\rangle
)=b(w_{i,a}\langle \sigma _{\alpha (i)}\rangle
w_{i,a}^{-1})=b(w_{i,a})w_{\alpha (i),b}\langle \sigma _{\beta (\alpha
(i))}\rangle w_{\alpha (i),b}^{-1}b(w_{i,a})^{-1}$$
Hence the map of
$Aut(W)$ into the group of permutations of $(1,\ldots ,m)$ defined
by $a\mapsto \alpha $ is a homomorphism. If $K$ is the kernel of
this homomorphism, then $K$ has finite index in $Aut (W)$ and for
all $a\in K$, $a(\langle \sigma _{i}\rangle )=w_{i,a}\langle \sigma
_{i}\rangle w_{i,a}^{-1}$ for all $i$.

For each $a\in K$ and for each $i$, we see that  $a(x)=w_{i,a}\tau
_{i}(x)w_{i,a}^{-1}$ for all $x\in \langle \sigma _{i}\rangle$,
where $\tau _{i}\in Aut(\langle\sigma _{i}\rangle)$. The map
$q_{i}:K\to Aut (\langle\sigma _{i}\rangle)$ defined by $a\mapsto
\tau_{i}$ is a homomorphism. The main result of \cite{HRT} shows
that $I(\langle\sigma _{i}\rangle)$, the inner automorphism group
of $\langle \sigma _{i}\rangle $ has finite index in $Aut
(\langle\sigma _{i}\rangle)$. Hence $C\equiv \cap
_{i=1}^{m}q_{i}^{-1}(I(\langle\sigma _{i}\rangle))$ is a subgroup
of finite index in $Aut (W)$.
\end{proof}

\begin{example}
This example has a subgroup generated by an edge in the Coxeter
diagram that is conjugate to that of a disjoint edge.
$$
\langle a,b,c,d: a^2=\cdots =d^2=1, (ab)^2, (ac)^2, (cd)^2,
(bc)^3, (ad)^3\rangle
$$
The element $bc$ conjugates $\langle a,b\rangle$ to $\langle
a,c\rangle$ and the element $ad$ conjugates $\langle a,c\rangle$
to $\langle c,d\rangle$. Hence the subgroup $\langle a,b\rangle$
is conjugate to $\langle c,d\rangle$.  Neither group is a maximal
FA subgroup.
\end{example}

\section{Visually Stable Subgroups}\label{stable}

If $(W,S)$ is a finitely generated Coxeter system, and $A$ is a
special subgroup for this system, then $A$ is {\it W-visually
stable} (or $W$-$VS$) if for any other Coxeter system $(W,S')$ for
$W$, $A$ is  conjugate to a special subgroup for $(W,S')$.

 Knowledge of the visually stable subgroups of a
Coxeter group is of interest in several important questions
related to the isomorphism problem for Coxeter groups. In
particular, this knowledge is useful for ``rigidity'' questions
(see for example \cite{charneydavis} and \cite{BMMN} and the
references there) and questions about when reflections are
preserved when passing between different Coxeter systems for $W$
(see \cite{BM}).

Maximal finite special subgroups are visually stable (see
\cite{Franzen}). Clearly $W$ and the trivial group are $W$-$VS$.
The following elementary observation is useful.

\begin{lemma}\label{tower} If $H$ is $W$-$VS$ and $K$ is $H$-$VS$ then $K$ is
$W$-$VS$.
\end{lemma}

The next result is a direct corollary to theorem
\ref{maxfatheorem}.

\begin{corollary} The maximal FA subgroups of a finitely generated
Coxeter group $W$ are $W$-$VS$.
\end{corollary}

Suppose $(W,S)$ is finitely generated. By Dunwoody's accessibility
theorem, $W$ decomposes as a graph of groups $\Lambda$, with
finite edge groups such that each vertex group is finite or
1-ended. Suppose $\Psi$ is a reduced visual decomposition for
$(W,S)$ derived from $\Lambda$ as given in theorem
\ref{maintheorem}. Any finite or 1-ended subgroup of $W$ is a
subgroup of a conjugate of a vertex group of $\Psi$ and so by
theorem \ref{closedunwoody} there is a bijection $\gamma$ from
the vertices of $\Lambda$ to the vertices of $\Psi$ such that
$\Lambda (V)$ is conjugate to $\Psi (\gamma (V))$. Note that the
vertex groups of $\Psi$ are finite or maximal 1-ended special
subgroups of $(W,S)$. In particular we have

\begin{proposition}\label{maxone}
The maximal 1-ended special subgroups of $(W,S)$ are $W$-$VS$.
\end{proposition}

The results of \cite{MTACCESS} and \cite{MTJSJ} imply that vertex
groups of the JSJ decompositions and strong accessibility
splittings considered there are $VS$. One example of these
implications is the following

\begin{proposition} \label{vsjsj} If $\Psi$ is a visual and irreducible
with respect to 2-ended splittings decomposition for a 1-ended
Coxeter system $(W,S)$ and $V$ is a vertex of $\Psi$, then
$\Psi(V)$ is $W$-$VS$.
\end{proposition}

An interesting situation occurs during repeated applications of
Lemma \ref{tower}, Proposition \ref{maxone} and Proposition
\ref{vsjsj} to an arbitrary finitely generated Coxeter system.
Suppose $H_1$ is a maximal 1-ended  special subgroup for $(W,S)$,
and $H_2$ is a maximal 1-ended special subgroup of a vertex group
of $\Psi$, a visual irreducible with respect to 2-ended
splittings decomposition of $H_1$. Then $H_2$ may split over a
2-ended group (just not in a way compatible with $\Psi$). The
vertex groups of an irreducible with respect to 2-ended
splittings decomposition of $H_2$ are $W$-$VS$. Continuing on
this line we have

\begin{theorem}
Suppose $(W,S)$ is a finitely generated Coxeter system, and
$A\subset S$ is maximal in the set of all $A'\subset S$ such that
the induced diagram for $A'$ is not separated by a subdiagram for
a finite or 2-ended special subgroup for $(W,S)$. Then $\langle
A\rangle$ is $W$-$VS$.
\end{theorem}

\section{A Final Application}\label{applicationssection}

Theorem \ref{maintheorem} and Corollary \ref{corends} can be
applied to ``visually'' characterize virtually free Coxeter
groups.

\begin{theorem}
The following are equivalent for any Coxeter system $(W,S)$ with
presentation diagram $\Gamma(W,S)$:
\begin{enumerate}
\item $W$ is virtually free, \item $W$ has a visual graph of
groups decomposition in which each vertex group is finite, \item
a) every complete subgraph of $\Gamma(W,S)$ is that of a finite
Coxeter group, and b) no induced subgraph of $\Gamma(W,S)$ is a
circuit of more than three vertices.
\end{enumerate}
\end{theorem}

\medskip
\noindent{\bf Proof:} Virtually free groups cannot contain 1-ended
subgroups. So by Dunwoody's accessibility theorem, every finitely
generated virtually free group is the fundamental group of a graph
of groups with finite vertex groups. By theorem \ref{maintheorem},
(1) implies (2). Conversely, any graph of groups with finite
vertex groups has fundamental group which is virtually free, and
(2) implies (1).

Again, as $W$ contains no 1-ended subgroup, corollary
\ref{corends} can be applied to show (2) implies (3a).  A circuit
of more than three vertices determines a 1-ended group.

Now suppose condition (3) holds in $\Gamma(W,S)$.  If
$\Gamma(W,S)$ were a complete graph then by a) it would be finite.
Instead take $x$ and $y$ in $\Gamma(W,S)$ not connected by an
edge. Consider a component $K$ of the complement of $st(x)$ in
$\Gamma(W,S)$. Let $S_{1}$ be the vertices in $st(x)$ adjacent to
a vertex in $K$. We claim that any two vertices $a$ and $b$ in
$S_{1}$ are connected by an edge in $\Gamma(W,S)$, otherwise a
minimal length path from $a$ to $b$ in $S_{1}\cup K$ meeting
$S_{1}$ only at its endpoints, together with the edges from $b$ to
$x$ and $x$ to $a$ would be a circuit in $\Gamma(W,S)$ of more
than three vertices contradicting b). Thus $S_{1}$ is a complete
graph, and so generates a finite subgroup $C$. Since $S_{1}$
separates $\Gamma(W,S)$, $W$ splits as a nontrivial visual
amalgamated product $A*_{C}B$. Inductively, (3) implies (2).
$\square$

\end{document}